\documentclass[a4paper,12pt]{amsart}
\usepackage{amscd}

\def\Cal{\mathcal}
\def\LOT{{\rm LOT}}
\newcommand{\nn}[1]{(\ref{#1})}
\newcommand{\ce}{{\Cal E}}
\newtheorem{theorem}{Theorem}[section]

\usepackage{amsmath,amsfonts,amssymb}

\def\sideremark#1{\ifvmode\leavevmode\fi\vadjust{\vbox to0pt{\vss
 \hbox to 0pt{\hskip\hsize\hskip1em
 \vbox{\hsize3cm\tiny\raggedright\pretolerance10000
 \noindent #1\hfill}\hss}\vbox to8pt{\vfil}\vss}}}%

\begin{document}
\renewcommand{\today}{}
\title{The conformal deformation detour complex for the obstruction tensor}
\author{Thomas P.\ Branson and A.\ Rod Gover}
\footnote{The second named author dedicates the paper to the memory of
Thomas P. Branson (1953 - 2006)}

\address{Department of Mathematics, The 
University of Iowa, Iowa City IA
52242 USA}

\address{Department of Mathematics\\
  The University of Auckland\\
  Private Bag 92019\\
  Auckland 1\\
  New Zealand} \email{gover@math.auckland.ac.nz}

\vspace{10pt}

\renewcommand{\arraystretch}{1}
\renewcommand{\arraystretch}{1.5}

\pagestyle{myheadings}

\begin{abstract}
  On pseudo-Riemannian manifolds of even dimension $n\geq 4$, with
  everywhere vanishing (Fefferman-Graham) obstruction tensor, we
  construct a complex of conformally invariant differential operators.
  The complex controls the infinitesimal deformations of
  obstruction-flat structures, and, in the case of Riemannian
  signature the complex is elliptic.
\end{abstract}
\maketitle \markboth{Branson \& Gover}{A deformation detour complex}
 \section{Introduction} 

 Ever since its introduction, as the natural tensor controlling a
conformal {\em Weyl relativity} \cite{Bach}, the Bach tensor has been
considered an intriguing and important natural invariant in
4-dimensional Riemannian and pseudo-Riemannian geometry.  In higher
even dimensions $n$, an analogue of the Bach tensor was discovered by
Fefferman and Graham \cite{FGast}; it arose as an obstruction to their
ambient metric construction. This {\em Fefferman-Graham obstruction
tensor}, which we denote $\mathcal{B}_{ab}$, is a natural tensor that
shares many of the properties of the Bach tensor.  It is a conformally
invariant, trace-free and divergence-free symmetric 2-tensor that
vanishes for conformally Einstein metrics.  Recently Graham and
Hirachi \cite{GrH} have shown that, also generalising the situation in
dimension 4, $\mathcal{B}_{ab}$ is the total metric variation of $\int
Q$, where $Q$ is Branson's Q-curvature \cite{tomsharp}. Along the way
they develop many of the properties of $\mathcal{B}_{ab}$. These
properties are also developed in \cite{GrH} via an alternative, but
equivalent, definition of the obstruction; it is given as the image of
a certain conformal operator when applied to the Weyl curvature.

 On conformally flat manifolds a large class of conformally invariant
 operators between irreducible bundles are organised into differential
 complexes known as generalised Bernstein-Gelfand-Gelfand (gBGG)
 complexes.  Recently there has been spectacular progress in the
 construction of curved analogues of these and analogous operators for
 other parabolic geometries, see \cite{CSSannals,CD} and references
 therein.  The issue of when these or other invariant operators may be
 combined to yield a complex, in some curved setting, is more subtle
 and there are only limited results.  On half-flat conformal
 4-manifolds there is a well-known systematic construction of a class
 of complexes. This has been extended to other structures by various
 authors and, recently, a rather general theory produced by \v Cap and
 Sou\v cek \cite{CapSou}. Some time ago it was observed by the current
 authors that, from the classification of conformally flat operators
 on the sphere \cite{BoC,ESlo}, it follows easily that on conformally
 flat structures there is a huge class of complexes related to gBGG
 complexes. We term these detour complexes, and each such is closely
 related to a gBGG complex. In \cite{deRham} we show that for the de
 Rham complex the corresponding detours yield conformally invariant
 differential complexes even in the fully curved case. These complexes
 are related to a conformal Hodge theory and yield new global
 invariants \cite{pontr}.  In \cite{GoPetob}, as background to the
 treatment of the obstruction tensor, some examples of other detour
 complexes are discussed in the conformally flat setting. These are
 related to deformations of conformally flat structures. It is in
 general a severe condition to require that a sequence of differential
 operators be both a complex, and conformally invariant; in general we
 do not expect complexes from conformally flat structures to admit
 curved analogues.  Nevertheless we show that one of the complexes
 from \cite{GoPetob} may be extended to curved structures that satisfy
 the curvature condition of everywhere vanishing Fefferman-Graham
 obstruction tensor. Background for much of the discussion here may be
 found in \cite{Tsrni}. Discussions with Andrew Waldron on related
 issues have been much appreciated.

\section{The deformation detour complex for the obstruction}

We work on a pseudo-Riemannian manifold $M$ of even dimension $n\geq
4$.  The metric will be denoted by $g$.  All bundles, sections spaces
and structures will be smooth. Let us write $T$ and $T^*$ for,
respectively, the tangent bundle and the cotangent bundle (although we
will also identify these via the metric without further comment). Then
$S^2_0T^*$ denotes the trace-free part of the second symmetric tensor
power of $T^*$. We will use the same notation for section spaces.

On a tangent vector field $X$, the Killing operator is the first order
differential operator $X\mapsto \mathcal{L}_X g$, where
$\mathcal{L}_X$ denotes the Lie derivative along the flow of $X$. We write 
$K X$ as a shorthand for this Killing operator on $X$ and, 
taking traces via the metric and its inverse, we write $K_0 X$ for the
trace-free part of $K X$. The operator $ X\mapsto K_0 X$ is known as
the conformal Killing operator and it is well-known that this is
conformally invariant. We write $K_0^*$ for the formal adjoint of this
operator.

Let us write $B$ (or sometimes $B^g$) for the linearisation, at the
metric $g$, of the operator which takes metrics $g$ to
$\mathcal{B}^g$. Since $\mathcal{B}$ is trace-free symmetric it follows that 
$B$ also takes values in $S^2_0T^*$. 
The main result is the following. 
\begin{theorem}\label{obflat} On pseudo-Riemannian manifolds with the 
Fefferman-Graham obstruction tensor vanishing everywhere,  the sequence 
of operators 
\begin{equation}\label{defdet}
 T\stackrel{K_0}{\to} S^2_0T^*\stackrel{B}{\to}
S^2_0T^* \stackrel{K_0^*}{\to} T
\end{equation}
is a formally self-adjoint complex of conformally invariant
operators. In Riemannian signature the complex is elliptic. 
\end{theorem}
\noindent{\bf Proof:} 
Since the conformal Killing
operator $K_0$ is conformally invariant, so is its formal
adjoint $K_0^*$. 
Viewing $\mathcal{B}$ as a tensor for each
metric (rather than a density-valued tensor) we have
\begin{equation}\label{Btrans}
\mathcal{B}^{\widehat{g}} = e^{(2-n)\Upsilon} \mathcal{B}^{g}
\end{equation}
where $\widehat{g}=e^{2\Upsilon}g$ for a smooth function $\Upsilon$. Thus 
setting $\Upsilon =t \omega$, for a real parameter $t$, and 
differentiating with respect to $t$ at $t=0$, we obtain the result that,
$$
B(\omega g)=\frac{2-n}{2}\omega \mathcal{B}^g
$$ for any function $\omega$. Thus if $\mathcal{B}=0$ then $B(\omega g)=0$;
since $B$ annihilates pure trace terms it descends to an operator
$S^2_0T^*{\to} S^2_0T^*$. From \nn{Btrans}, $B^g$ is conformally
bi-invariant as follows,
\begin{equation}\label{bi}
B^{\widehat{g}} (e^{2\Upsilon} h)= e^{(2-n)\Upsilon} B^{g} h ~,
\end{equation}
for $h$ a section in $S^2_0T^*$.
 So the sequence is conformally invariant as claimed. 

Since the obstruction tensor $\mathcal{B}^g$ is 
a natural Riemannian invariant, the Lie
derivative along a vector field $X$ is 
$$
\mathcal{L}_X \mathcal{B}= B(\mathcal{L}_X g)=  B (K X),
$$ by the chain rule. We have seen that, under the curvature condition 
$\mathcal{B}=0$, $B$ annihilates  terms that are 
pure trace and so we have $B(K X)=B(K_0 X)=0$.
That is,  the composition $B K_0$ is zero. 

As mentioned above, the obstruction tensor
$\mathcal{B}$ arises as the total metric variation of the integral of
$Q$. Thus its linearisation is a second variation of $\int Q$ and
hence is a formally self-adjoint operator on the space of metric
variations.  So $B$ is formally self-adjoint on $S^2_0 T^*$. Thus the
sequence \nn{defdet} is formally self-adjoint and hence a complex.

We write $\nabla$ for the Levi-Civita connection, and $\Delta$ for the
Laplacian $g^{ab}\nabla_a\nabla_b$ (given here in abstract index
notation).  Note that from \cite{GrH,GoPetob} we have
$$
\mathcal{B}=\Delta^{n/2-2}\delta^\nabla \delta^\nabla C + \LOT ~.
$$ 
where $C$ is the Weyl curvature tensor, $\delta^\nabla$ indicates the obvious 
covariant divergence operators (in an abstract index notation $\delta^\nabla \delta^\nabla C$ is $\nabla^b\nabla^d C_{abcd}$),  and $\LOT$ indicates lower order terms. Thus at leading order its linearisation takes the form 
$$
B= \Delta^{n/2-2}\delta^\nabla \delta^\nabla W + \LOT ~.
$$ where $W$ is the linearisation of the Weyl curvature, as a
 differential operator on $S^2_0 T^*$. But from elementary classical
 invariant theory it is easily verified that the double divergence is
 the formal adjoint of the operator $W$. This may be commuted past the
 Laplacian power at the expense of lower order terms and so the
 leading term of $B$ is simply
$$
W^* \Delta^{n/2-2} W
$$ and whence, in the case of Riemannian signature, the leading symbol
of $B$, $\sigma(B)$ has the same kernel as $\sigma (W)$. This kernel
is the same as the image space of the leading symbol of $K_0$ since these give
the leading terms of a resolution \cite{GG,CD,GoPetob}
$$ T\stackrel{K_0}{\to} S^2_0T^*\stackrel{W}{\to}
\mathcal{W} \to \cdots $$ 
that controls deformations of conformally
flat structures. Here $\mathcal{W}$ is the space of algebraic Weyl tensors.
 Thus the symbol sequence is exact and the complex
\nn{defdet} is elliptic as claimed.  \quad $\Box$

\noindent{\bf Remarks:} Since the operator $K_0$ takes 
infinitesimal diffeomorphisms to their action on
conformal structure, it is  
clear
that the first cohomology of the complex \nn{defdet} is the formal
tangent space to the moduli space of deformations of obstruction-flat
geometries. In fact the complex controls the full formal
deformation theory, but this will be taken up elsewhere. 

To simplify the exposition here we avoided the use of conformal
densities and discussion of conformal structures. Writing $E[w]$ for
the bundle of conformal densities of weight $w$ (see e.g.\ \cite{CapGoamb}), 
and 
$T^*[w]$ as a shorthand for $T^*\otimes E[w]$ and so forth, the complex is 
$$
T\stackrel{K_0}{\to} S^2_0T^*[2]\stackrel{B}{\to}
S^2_0T^*[2-n] \stackrel{K_0^*}{\to} T^*[-n]
$$ 
or, in the notation of \cite{GoPetob} (where the conformally flat version of this is discussed),
$$
\ce^1[2]\stackrel{K_0}{\to} \ce^{1,1}[2]\stackrel{B}{\to}
\ce_{1,1}[-2] \stackrel{K_0^*}{\to} \ce_{1}[-2]~.
$$

The main idea for the dimension 4 version of this complex was
developed during the Spring 2001 programme ``Spectral Invariants --
Analytical and Geometric Aspects'' at the Mathematical Sciences
Research Institute, Berkeley. This support by the MSRI is greatly
appreciated.  The second author would also like to thank the Royal
Society of New Zealand for support via Marsden Grant no.\ 02-UOA-108.

\end{document}